\documentclass{article}
\usepackage{amssymb}
\usepackage{amsbsy}
\begin{document}
\centerline{\Large \bf } \vskip 6pt

\begin{center}{\Large \bf Orlicz version of the mixed width integrals}\footnote{Research is supported by
National Natural Sciences Foundation of China (10971205, 11371334).}\end{center}
\vskip 6pt

\centerline{Chang-Jian Zhao}

\begin{center}{\it Department of Mathematics,
China Jiliang University, Hangzhou 310018, P. R.
China}\end{center} \centerline{\it Email: chjzhao@163.com}

\vskip 6pt

\begin{center}
\begin{minipage}{12cm}
{\bf Abstract}~ In the paper, our main aim is to generalize the
width integrals to the Orlicz space. Under the framework of Orlicz
Brunn-Minkowski theory, we introduce a new affine geometric
quantity by calculating Orlicz first order variation of the width
integrals, and call as Orlicz mixed width integrals. The
fundamental notions and conclusions of the width integrals and
Minkoswki and Brunn-Minkowski inequalities for the width integrals
are extended to an Orlicz setting and the related concepts and
inequalities of $L_{p}$-mixed width integrals of convex body are
also derived.


\end{minipage}
\end{center}
\vskip 20pt

\noindent{\large \bf 1 ~Introducation}\vskip 10pt

One of the most important operations in geometry is vector
addition. As an operation between sets $K$ and $L$, defined by
$$K+L=\{x+y: x\in K,y\in L\},$$ it is usually called Minkowski
addition and plays an important role in the
Brunn-Minkowski theory. During the last few decades, the theory
has been extended to $L_{p}$-Brunn-Minkowski theory. For compact convex
sets $K$ and $L$ containing the origin in its interior and $1\leq p<\infty$, the $L_p$ addition of
$K$ and $L$, introduced by Firey in [3] or [4], is defined by
$$h(K+_{p}L,x)^{p}=h(K,x)^{p}+h(L,x)^{p},\eqno(1.1)$$
for all $x\in {\Bbb R}^{n}$ and compact convex sets $K$ and $L$ in
${\Bbb R}^{n}$ containing the origin. When $p=\infty$, (1.1) is
interpreted as $h(K+_{\infty}L,x)=\max\{h(K,x), h(L,x)\}$. Here the functions are the support functions. If $K$ is
a nonempty closed (not necessarily bounded) convex set in ${\Bbb
R}^{n}$, then the support function $h(K, x)$ of $K$ is defined by
$$h(K,x)=\max\{x\cdot y: y\in K\},\eqno(1.2)$$
for $x\in {\Bbb R}^{n}$. A nonempty closed convex set is uniquely determined by its
support function. $L_{p}$ addition and inequalities are the
fundamental and core content in the $L_{p}$ Brunn-Minkowski
theory. For recent important results and more information from
this theory, refer to [8], [9], [10], [11], [15], [16], [18], [19],
[20], [21], [22], [25], [26], [29], [30] and [31]. In recent years, a new extension of
$L_{p}$-Brunn-Minkowski theory is to Orlicz-Brunn-Minkowski
theory, initiated by Lutwak, Yang, and Zhang [23] and [24].
Gardner, Hug and Weil [6] constructed a general framework for the
Orlicz-Brunn-Minkowski theory, and made clear for the first time
the relation to Orlicz spaces and norms. The Orlicz addition of
convex bodies was introduced and the Orlicz Brunn-Minkowski inequality is obtained (see [32]).
The Orlicz centroid inequality for convex bodies was introduced in
[39] which is an extension from star to convex bodies. The
Orlicz-Brunn-Minkowski theory and its dual theory have attracted
people's attention. The other articles to promote the theory
can be found in literatures [7], [13], [14], [27], [33], [38],
[34], [35], [36] and [37].

For $u\in S^{n-1}$, the half width of convex body $K$ in the
direction $u$, defined by.
$$b(K,u)=\frac{1}{2}(h(K,u)+h(K,-u)).\eqno(1.3)$$
Convex bodies $K,L$ are said to have similar width if there exist a
constant $\lambda>0$ such that $b(K,u)=\lambda b(L,u),$ for all
$u\in S^{n-1}$. Lutwak [17] introduced the width integrals
: For
$0\leq i<n$, the width integral of convex body $K$, denotes by
$A_{i}(K),$ defined by
$$A_{i}(K)=\frac{1}{n}\int_{S^{n-1}}b(K,u)^{n-i}dS(u).\eqno(1.4)$$

In the paper, our main aim is to generalize the width integrals to
Orlicz space. Under the framework of Orlicz Brunn-Minkowski
theory, we introduce a new affine geometric quantity-Orlicz mixed width integrals. The fundamental notions and
conclusions of the width integrals and Minkoswki and
Brunn-Minkowski inequalities for the width integrals are extended
to an Orlicz setting. The related concepts and inequalities of
$L_{p}$-mixed width integrals of convex body are also derived.
The new Orlicz Minkowski and Brunn-Minkowski inequalities in
special case yield the $L_p$-dual Minkowski, and Brunn-Minkowski
inequalities for the $L_{p}$-mixed width integrals.

In Section 3, we introduce a notion of Orlicz width addition
$K+_{\phi}L$ of convex bodies $K$ and $L$, defined by
$$\phi\left(\frac{b(K,x)}{b(K+_{\phi}L,x)}
,\frac{b(L,x)}{b(K+_{\phi}L,x)}\right)=1.\eqno(1.5)$$ Here
$\phi\in \Phi_{2}$, the set of convex function
$\phi:[0,\infty)^{2}\rightarrow (0,\infty)$ that are decreasing in
each variable and satisfy $\phi(0,0)=\infty$ and
$\phi(\infty,1)=\phi(1,\infty)=1$. The particular instance is
$\phi(x_{1},x_{2})=\phi_{1}(x_{1})+\varepsilon\phi_{2}(x_{2})$ for
$\varepsilon>0$ and some $\phi_{1},\phi_{2}\in\Phi$, where the
sets of convex functions $\phi_{1}, \phi_{2}:[0,\infty)\rightarrow
(0,\infty)$ that are decreasing and satisfy
$\phi_{1}(0)=\phi_{2}(0)=\infty$,
$\phi_{1}(\infty)=\phi_{2}(\infty)=0$ and
$\phi_{1}(1)=\phi_{2}(1)=1$ .

Complying with the basic spirit of Aleksandrov
[1], mixed quermassintegrals [2] and $L_p$-mixed quermassintegrals [16], we concentrate on the study of the first order Orlicz variational of the
width integrals. In Section 4, we prove that the first order
Orlicz variation of the width integrals can be expressed as: For
convex bodies $K,L$, and $\phi_{1},\phi_{2}\in {\Phi}$, $0\leq
i<n$ and $\varepsilon>0$,
$$\frac{d}{d\varepsilon}\bigg|_{\varepsilon=0^{+}}A_{i}
(K+_{\phi}\varepsilon\cdot L)=\frac{n-i}{(\phi_{1}) '_{r}(1)}\cdot
A_{\phi_{2},i}(K,L).\eqno(1.6)$$ In this first order variational
equation (1.6), we define a new geometric quantity,
Orlicz mixed width integrals $A_{\phi_{2},i}(K,L)$, is defined by
$$A_{\phi_{2},i}(K,L):=\frac{(\phi_{1})
'_{r}(1)}{n-i}\cdot\frac{d}{d\varepsilon}\bigg|_{\varepsilon=0^{+}}A_{i}
(K+_{\phi}\varepsilon\cdot L).\eqno(1.7)$$ We also prove the new
affine geometric quantity has an integral representation.
$$A_{\phi,i}(K,L)=\frac{1}{n}\int_{S^{n-1}}\phi
\left(\frac{b(L,u)}{b(K,u)}\right)b(K,u)^{n-i}dS(u).\eqno(1.8)$$
In Section 5, we establish an Orlicz Minkowski inequality for the
Orlicz mixed width integrals: If $K$ and $L$ are convex bodies,
$0\leq i<n$ and $\phi\in \Phi$, then
$$A_{\phi,i}(K,L)\geq
A_{i}(K)\cdot\phi\left(\left(\frac{A_{i}(L)}{A_{i}(K)}\right)^{1/(n-i)}\right)
.\eqno(1.9)$$ If $\phi$ is strictly convex, equality holds if and
only if $K$ and $L$ have similar width. In Section 6, we establish
an Orlicz Brunn-Minkowski inequality for the Orlicz width addition
and the width integrals. If $K$ and $L$ are convex bodies, $0\leq
i<n$ and $\phi\in \Phi_{2}$, then
$$1\geq\phi\left(\left(\frac{A_{i}(K)}
{A_{i}(K+_{\phi}L)}\right)^{1/(n-i)},\left(\frac{A_{i}(L)}{A_{i}(K+_
{\phi}L)}\right)^{1/(n-i)}\right).\eqno(1.10)$$ If $\phi$ is
strictly convex, equality holds if and only if $K$ and $L$ have
similar width.

\vskip 10pt \noindent{\large \bf 2 ~Preliminaries}\vskip 10pt

The setting for this paper is $n$-dimensional Euclidean space
${\Bbb R}^{n}$. A body in ${\Bbb R}^{n}$ is a compact set equal to
the closure of its interior. A set $K$ is called a convex body if
it is compact and convex subsets with non-empty interiors. Let
${\cal K}^{n}$ denote the class of convex bodies containing the
origin in their interiors in ${\Bbb R}^{n}$. For a compact set
$K\subset {\Bbb R}^{n}$, we write $V(K)$ for the ($n$-dimensional)
Lebesgue measure of $K$ and call this the volume of $K$.  The
letter $B$ denotes the unit ball centered at the origin. The
support function of convex body $K$ is homogeneous of degree $1$,
that is (see e.g. [28]),
$$h(K,ru)=rh(K,u),$$
for all $u\in {S}^{n-1}$ and $r>0$. Let $\delta$ denote the
Hausdorff metric, as follows, if $K, L\in {\cal K}^{n}$, then
$$\delta(K,L)=|h(K,u)-h(L,u)|_{\infty}.$$ If $K\in {\cal K}^{n}$ and
$A\in{\rm GL(n)}$, then (see
e.g. [5], p.17)
$$h(AK,x)=h(K, A^{t}x),\eqno(2.1)$$
for all $x\in {\Bbb R}^{n}$.

For $K_{i}\in {\cal K}^{n}, i=1,\ldots,m$, define the real numbers
$R_{K_{i}}$ and $r_{K_{i}}$ by
$$R_{K_{i}}=\max_{u\in S^{n-1}}b(K_{i},u),~~ {\rm and}~~ r_{K_{i}}=\min_{u\in S^{n-1}}b(K_{i},u).\eqno(2.2)$$
Obviously, $0<r_{K_{i}}<R_{K_{i}},$ for all $K_{i}\in {\cal
S}^{n}$, and writing $R=\max\{R_{K_{i}}\}$ and
$r=\min\{r_{K_{i}}\}$, where $i=1,\ldots,m.$

\vskip 8pt {\it 2.1~ Mixed width integrals}\vskip 8pt

If $K_{1},\ldots,K_{n}\in {\cal K}^{n}$, the mixed width integral
of $K_{1},\ldots,K_{n}$, denoted by $A(K_{1},\ldots,K_{n})$,
defined by (see [17])
$$A(K_{1},\ldots,K_{n})=\frac{1}{n}\int_{S^{n-1}}b(K_{1},u)\cdots b(K_{n},u)dS(u).$$
If $K_{1}=\cdots=K_{n-i}=K,$ $K_{n-i+1}=\cdots=K_{n}=B$, the mixed
width integral $A(K_{1},\ldots,K_{n})$ is written as
$A_{i}(K)$
and call width integral of $K$. Obviously, For $K\in {\cal K}^{n}$
and $0\leq i<n$, we have
$$A_{i}(K)=\frac{1}{n}\int_{S^{n-1}}b(K,u)^{n-i}dS(u).\eqno(2.3)$$
If $K_{1}=\cdots=K_{n-i-1}=K,$ $K_{n-i}=\cdots=K_{n-1}=B$ and
$K_{n}=L$, the mixed width integral
$A(\underbrace{K,\ldots,K}_{n-i-1},\underbrace{B,\ldots,B}_{i},L)$
is written as $A_{i}(K,L)$ and call $i$-th mixed width integral of
$K$ and $L$. For $K,L\in {\cal K}^{n}$ and $0\leq i<n$, it is easy
that
$$A_{i}(K,L)=\frac{1}{n}\int_{S^{n-1}}b(K,u)^{n-i-1}b(L,u)dS(u).\eqno(2.4)$$
This integral representation (2.4), together with H\"{o}lder
inequality, immediately gives: The Minkowski inequality for the
$i$-th mixed width integral. If $K,L\in {\cal K}^{n}$ and $0\leq
i<n$, then
$$A_{i}(K,L)^{n-i}\leq A_{i}(K)^{n-i-1}A_{i}(L),\eqno(2.5)$$
with equality if and only if $K$ and $L$ have similar width.

\vskip 8pt {\it 2.2~ $L_p$-mixed width integrals}\vskip 8pt

Putting $\phi(x_{1},x_{2})=x_{1}^{-p}+x_{2}^{-p}$ and $p\geq 1$ in
(1.5), the Orlicz width addition $+_{\phi}$ becomes the $L_p$-width addition, denoted by $+_{p}$, and call as $L_p$-width
addition of convex bodies $K$ and $L$.
$$b(K+_{p}L,u)^{-p}=b(K,u)^{-p}+b(L,u)^{-p},\eqno(2.6)$$
for $u\in S^{n-1}$. The following result follows immediately form
(2.6) with $p\geq 1$.
$$-\frac{np}{n-i}\lim_{\varepsilon\rightarrow 0^{+}}\frac{A_{i}(K+_{p}
\varepsilon\cdot
L)-A_{i}(L)}{\varepsilon}=\frac{1}{n}\int_{S^{n-1}}b(K,u)^{n-i+p}b(L,u)^{-p}dS(u).$$

{\bf Definition 2.1}~ Let $K,L\in {\cal K}^{n}$, $0\leq i<n$ and
$p\geq 1$, the $L_p$-width integral of convex bodies $K$ and $L$,
denoted by $A_{-p,i}(K,L)$, defined by
$$A_{-p,i}(K,L)=\frac{1}{n}\int_{S^{n-1}}b(K,u)^{n-i+p}b(L,u)^{-p}dS(u).\eqno(2.7)$$
Obviously, when $K=L$, the $L_p$-mixed width integral
$A_{-p,i}(K,K)$ becomes the width integral $A_{i}(K).$ This
integral representation (2.7), together with H\"{o}lder
inequality, immediately gives:

{\bf Proposition 2.2}~ {\it If $K,L\in {\cal K}^{n}$, $0\leq i<n$
and $p\geq 1$, then
$$A_{-p,i}(K,L)^{n-i}\geq A_{i}(K)^{n-i+p}A_{i}(L)^{-p},\eqno(2.8)$$
with equality if and only if $K$ and $L$ have similar width.}

{\bf Proposition 2.3}~ {\it If $K,L\in {\cal K}^{n}$, $0\leq i<n$
and $p\geq 1$, then
$$A_{i}(K+_{p}L)^{-p/(n-i)}\geq A_{i}(K)^{-p/(n-i)}+A_{i}(L)^{-p/(n-i)},\eqno(2.9)$$
with equality if and only if $K$ and $L$ are dilates.}

{\it Proof}~ From (2.6) and (2.7), it is easily seen that the
$L_p$-width integrals is linear with respect to the $L_p$-width
addition, and together with inequality (2.8) show that for $p\geq
1$
$$A_{-p,i}(Q, K+_{p}L)=A_{-p,i}(Q,K)+A_{-p,i}(Q,L)~~~~~~~~~~~~~~~~~~~~~~~~~~~~~~~~~~~~~~~~~~~~$$
$$~~~~~~~~~~~~
~\geq A_{i}(Q)^{(n-i+p)/(n-i)}(A_{i}(K)^{-p/(n-i)}+A_{i}(L)^{-p/(n-i)}),$$
with equality if and only if $K$ and $L$ have similar width.

Take $K+_{p}L$ for $Q$, recall that $A_{-p,i}(Q,Q)=A_{i}(Q)$,
inequality (2.9) follows
easy.\hfill$\Box$

\vskip 10pt \noindent{\large \bf 3 ~Orlicz width addition}\vskip
10pt

Throughout the paper, the standard orthonormal basis for ${\Bbb
R}^{n}$ will be $\{e_{1},\ldots,e_{n}\}$. Let $\Phi_{n},$
$n\in{\Bbb N}$, denote the set of convex function
$\phi:[0,\infty)^{n} \rightarrow (0,\infty)$ that are strictly
decreasing in each variable and satisfy $\phi(0)=\infty$ and
$\phi(e_{j})=1$, $j=1,\ldots,n$. When $n=1$, we shall write $\Phi$
instead of $\Phi_{1}$. The left derivative and right derivative of
a real-valued function $f$ are denoted by $(f)'_{l}$ and
$(f)'_{r}$, respectively. We first define the Orlicz width
addition.

{\bf Definition 3.1}~ Let $m\geq 2, \phi\in\Phi_{m}$, $K_{j}\in
{\cal K}^{n}$ and $j=1,\ldots,m$, define the Orlicz width addition
of $K_{1},\ldots,K_{m}$, denotes by
$+_{\phi}(K_{1},\ldots,K_{m})$, defined by
$$b(+_{\phi}(K_{1},\ldots,K_{m}),u)=\sup\left\{\lambda>0:
\phi\left(\frac{b(K_{1},u)}{\lambda},\ldots,\frac{b(K_{m},u)}{\lambda}\right)\leq
1\right\},\eqno(3.1)$$ for $u\in S^{n-1}.$ Equivalently, the
Orlicz width addition $+_{\phi}(K_{1},\ldots,K_{m})$ can be
defined implicitly by
$$\phi\left(\frac{b(K_{1},u)}{b(+_{\phi}(K_{1},
\ldots,K_{m}),u)},\ldots,\frac{b(K_{m},u)}{b(+_{\phi}(K_{1},
\ldots,K_{m}),u)}\right)=1,\eqno(3.2)$$ for all $u\in {S}^{n-1}$.

An important special case is obtained when
$$\phi(x_{1},\ldots,x_{m})=\sum_{j=1}^{m}\phi_{j}(x_{j}),$$
for some fixed $\phi_{j}\in \Phi$ such that
$\phi_{1}(1)=\cdots=\phi_{m}(1)=1$. We then write
$+_{\phi}(K_{1},\ldots,K_{m})=K_{1}+_{\phi}\cdots+_{\phi}K_{m}.$
This means that $K_{1}+_{\phi}\cdots+_{\phi}K_{m}$ is defined
either by
$$b(K_{1}+_{\phi}\cdots+_{\phi}K_{m},u)=
\sup\left\{\lambda>0:\sum_{j=1}^{m}\phi_{j}\left(\frac{b(K_{j},u)}{\lambda}
\right)\leq 1\right\},\eqno(3.3)$$ for all $u\in {S}^{n-1}$, or by
the corresponding special case of (3.2).

{\bf Lemma 3.2}~ {\it The Orlicz width addition $+_{\phi}: ({\cal
K}^{n})^{m}\rightarrow {\cal K}^{n}$ is monotonic and has the
identity property.}

{\it Proof}~ Suppose $K_{j}\subset L_{j}$, $j=1,\ldots,m$, where
$K_{j}, L_{j}\in {\cal K}^{n}.$ By using (3.1), and in view of
$K_{1}\subset L_{1}$ and the fact that $\phi$ is decreasing in the
first variable, we obtain
$$b(+_{\phi}(L_{1},K_{2}
\ldots,K_{m}),u)~~~~~~~~~~~~~~~~~~~~~~~~~~~~~~~~~~~~~~~~~~~~~~~~~~~~~~~~~~~~~~~$$
$$=\sup\left\{\lambda>0:\phi\left(\frac{b(L_{1},u)}{\lambda},\frac{b(K_{2},u)}{\lambda},\ldots,\frac{b(K_{m},u)}{\lambda}\right)\leq
1\right\}~$$
$$\geq\sup\left\{\lambda>0:\phi\left(\frac{b(K_{1},u)}{\lambda},\frac{b(K_{2},u)}{\lambda},\ldots,\frac{b(K_{m},u)}{\lambda}\right)\leq 1\right\}$$
$$=b(+_{\phi}(K_{1},K_{2},
\ldots,K_{m}),u).~~~~~~~~~~~~~~~~~~~~~~~~~~~~~~~~~~~~~~~~$$ By
repeating this argument for each of the other $(m-1)$ variables,
we have $b(+_{\phi}(K_{1}, \ldots,K_{m}),u)\leq
b(+_{\phi}(L_{1},\ldots,L_{m}),u)$.

The identity property is obvious from
(3.2).\hfill$\Box$

{\bf Lemma 3.3}~ {\it The Orlicz width addition $+_{\phi}: ({\cal
K}^{n})^{m}\rightarrow {\cal K}^{n}$ is $GL(n)$ covariant.}

{\it Proof}~ From (2.1), (3.1) and let $A\in GL(n)$, we obtain
$$b(+_{\phi}(AK_{1},AK_{2}
\ldots,AK_{m}),u)~~~~~~~~~~~~~~~~~~~~~~~~~~~~~~~~~~~~~~~~~~~~~~~~~~~~~~~~$$
$$~~~~=\sup\left\{\lambda>0:\phi\left(\frac{b(AK_{1},u)}{\lambda},\frac{b(AK_{2},u)}{\lambda},\ldots,\frac{b(AK_{m},u)}{\lambda}\right)\leq
1\right\}$$
$$~~~~~~~=\sup\left\{\lambda>0:\phi\left(\frac{b(K_{1},A^{t}u)}{\lambda},\frac{b(K_{2},A^{t}u)}{\lambda},\ldots,\frac{b(K_{m},A^{t}u)}{\lambda}\right)\leq 1\right\}$$
$$=b(+_{\phi}(K_{1},
\ldots,K_{m}),A^{t}u)~~~~~~~~~~~~~~~~~~~~~~~~~~~~~~~~~~~~~~~~~~~~~$$
$$=b(+_{\phi}(K_{1},
\ldots,K_{m}),u).~~~~~~~~~~~~~~~~~~~~~~~~~~~~~~~~~~~~~~~~~~~~~~~.$$

This shows Orlicz width addition $+_{\phi}$ is $GL(n)$ covariant.
\hfill$\Box$

{\bf Lemma 3.4}~ {\it Suppose $K,\ldots,K_{m}\in{\cal K}^{n}$. If
$\phi\in \Phi^{m}$, then
$$\phi\left(\frac{b(K_{1},u)}{t}\right)+\cdots+\phi\left(\frac{b(K_{m},u)}{t}\right)=1$$
if and only if}
$$b(+_{\phi}(K_{1},
\ldots,K_{m}),u)=t$$

{\it Proof}~  This follows immediately from Definition
3.1.\hfill$\Box$

{\bf Lemma 3.5}~ {\it Suppose $K_{1},\ldots,K_{m}\in{\cal K}^{n}$.
If $\phi\in \Phi^{m}$, then}
$$\frac{r}{\phi^{-1}(\frac{1}{m})}\leq b(+_{\phi}(K_{1},
\ldots,K_{m}),u)\leq\frac{R}{\phi^{-1}(\frac{1}{m})}.$$

{\it Proof}~ Suppose $b(+_{\phi}(K_{1}, \ldots,K_{m}),u)=t$. From
Lemma 3.4 and noting that $\phi$ is strictly deceasing on $(0,
\infty)$, we have
$$1=\phi\left(\frac{b(K_{1},u)}{t}\right)+\cdots+\phi\left(\frac{b(K_{m},u)}{t}\right)$$
$$\leq\phi\left(\frac{r_{K_{1}}}{t}\right)+\cdots+\phi\left(\frac{r_{K_{m}}}{t}\right)~~~~~~~~~~~$$
$$\leq m\phi\left(\frac{r}{t}\right).~~~~~~~~~~~~~~~~~~~~~~~~~~~~~~~~~~$$
Noting that $\phi^{-1}$ is strictly deceasing on $(0,\infty)$, we
obtain the lower bound for $b(+_{\phi}(K_{1}, \ldots,K_{m}),u)$:
$$t\geq\frac{r}{\phi^{-1}(\frac{1}{m})}.$$
To obtain the upper estimate, observe that from the lemma 3.4,
together with the convexity and the fact $\phi$ is strictly
deceasing on $(0,\infty)$, we have
$$1=\phi\left(\frac{b(K_{1},u)}{t}\right)+\cdots+\phi\left(\frac{b(K_{m},u)}{t}\right)$$
$$\geq m\phi\left(\frac{b(K_{1},u)+\cdots+b(K_{m},u)}{mt}\right)~~$$
$$\geq m\phi\left(\frac{R}{t}\right).~~~~~~~~~~~~~~~~~~~~~~~~~~~~~~~~$$
Then we obtain the upper estimate:
$$t\leq\frac{R}{\phi^{-1}(\frac{1}{m})}.$$
\hfill$\Box$

{\bf Lemma 3.6}~ {\it The Orlicz width addition $+_{\phi}: ({\cal
K}^{n})^{m}\rightarrow {\cal K}^{n}$ is continuous.}

{\it Proof}~ To see this, indeed, let $K_{ij}\in {\cal K}^{n},$ $i
\in {\Bbb N}\cup\{0\},$ $j=1,\ldots,m,$ be such that
$K_{ij}\rightarrow K_{0j}$ as $i\rightarrow\infty$. Let
$$b(+_{\phi}
(K_{i1},\ldots,K_{im}),u)=t_{i}.$$ Then Lemma 3.5 shows
$$\frac{r_{ij}}{\phi^{-1}(\frac{1}{m})}\leq t_{i}\leq\frac{R_{ij}}{\phi^{-1}(\frac{1}{m})},$$
where $r_{ij}=\min\{r_{K_{ij}}\}$ and $R_{ij}=\max\{R_{K_{ij}}\}.$
Since $K_{ij}\rightarrow K_{0j}$, we have $R_{K_{ij}}\rightarrow
R_{K_{0j}}<\infty$ and $r_{K_{ij}}\rightarrow r_{K_{0j}}>0$, and
thus there exist $a,b$ such that $0<a\leq t_{i}\leq b<\infty$ for
all $i$. To show that the bounded sequence $\{t_{i}\}$ converges
to $b(+_{\phi}(K_{01},\ldots,K_{0m}),u),$ we show that every
convergent subsequence of $\{t_{i}\}$ converges to
$b(+_{\phi}(K_{01},\ldots,K_{0m}),u)$. Denote any subsequence of
$\{t_{i}\}$ by $\{t_{i}\}$ as well, and suppose that for this
subsequence, we have
$$t_{i}\rightarrow t_{*}.$$
Obviously $a\leq t_{*}\leq b.$ Noting that $\phi$ is continuous
function, we obtain
$$t_{*}\rightarrow\sup\left\{t_{*}>0:
\phi\left(\frac{b(K_{01},u)}{t_{*}},\ldots,\frac{b(K_{0m},u)}{t_{*}}\right)\leq
1\right\}$$
$$~~~~~~=b(+_{\phi}(K_{01},\ldots,K_{0m}),u).~~~~~~~~~~~~~~~~~~~~~~~~~~~~~~$$
Hence
$$b(+_{\phi}
(K_{i1},\ldots,K_{im}),u)\rightarrow
b(+_{\phi}(K_{01},\ldots,K_{0m}),u)$$ as $i\rightarrow\infty$.

This shows that the Orlicz width addition $+_{\phi}: ({\cal
K}^{n})^{m}\rightarrow {\cal K}^{n}$ is continuous.
\hfill$\Box$

Next, we define the Orlicz width linear combination on the case
$m=2$.

{\bf Definition 3.7}~ Orlicz width linear combination
$+_{\phi}(K,L,\alpha,\beta)$ for $K,L\in {\cal K}^{n}$, and
$\alpha,\beta\geq 0$ (not both zero), defined by
$$\alpha\cdot\phi_{1}\left(\frac{b(K,u)}{b(+_{\phi}(K,L,\alpha,\beta),u)}\right)+
\beta\cdot\phi_{2}\left(\frac{b(L,u)}
{b(+_{\phi}(K,L,\alpha,\beta),u)}\right)=1,\eqno(3.4)$$ for all
$u\in {S}^{n-1}$.

We shall write $K+_{\phi}\varepsilon\cdot L$ instead of
$+_{\phi}(K,L,1,\varepsilon)$, for $\varepsilon\geq 0$ and assume
throughout that this is defined by (3.1), if $\alpha=1,
\beta=\varepsilon$ and $\phi\in \Phi$. We shall write $K+_{\phi}L$
instead of $+_{\phi}(K,L,1,1)$ and call the Orlicz width addition
of $K$ and $L$.

\vskip 10pt \noindent{\large \bf 4 ~Orlicz mixed width
integrals}\vskip 10pt

In order to define Orlicz mixed width integrals, we need the
following lemmas.

{\bf Lemma 4.1}~ {\it Let $\phi\in \Phi$ and $\varepsilon>0$. If
$K,L\in {\cal K}^{n}$, then} $K+_{\phi}\varepsilon\cdot L\in {\cal
K}^{n}.$

{\it Proof}~ Let $u_{0}\in S^{n-1}$, for any subsequence
$\{u_{i}\}\subset S^{n-1}$ such that $u_{i}\rightarrow u_{0}.$ as
$i\rightarrow \infty.$

Let
$$b(K+_{\phi}
L,u_{i})=\lambda_{i}.$$ Then Lemma 3.5 shows
$$\frac{r}{\phi^{-1}(\frac{1}{2})}\leq \lambda_{i}\leq\frac{R}{\phi^{-1}(\frac{1}{2})},$$
where $R=\max\{R_{K},R_{L}\}$ and $r=\min\{r_{K},r_{L}\}.$

Since $K,L\in {\cal K}^{n}$, we have $0<r_{K}\leq R_{K} <\infty$
and $0<r_{L}\leq R_{L} <\infty$, and thus there exist $a,b$ such
that $0<a\leq \lambda_{i}\leq b<\infty$ for all $i$. To show that
the bounded sequence $\{\lambda_{i}\}$ converges to
$b(K+_{\phi}L,u_{0}),$ we show that every convergent subsequence
of $\{\lambda_{i}\}$ converges to $b(K+_{\phi}L,u_{0})$. Denote
any subsequence of $\{\lambda_{i}\}$ by $\{\lambda_{i}\}$ as well,
and suppose that for this subsequence, we have
$$\lambda_{i}\rightarrow \lambda_{0}.$$
Obviously $a\leq \lambda_{0}\leq b.$ From (3.4) and note that
$\phi_{1},\phi_{2}$ are continuous functions, so $\phi^{-1}_{1}$
is continuous, we obtain
$$\lambda_{i}\rightarrow\frac{b(K,u_{0})}{\displaystyle\phi_{1}^{-1}\left(1-\varepsilon\phi_{2}
\left(\frac{b(L,u_{0})}{\lambda_{0}}\right)\right)}$$ as
$i\rightarrow \infty.$ Hence
$$\phi_{1}\left(\frac{b(K,u_{0})}{\lambda_{0}}\right)
+\varepsilon\phi_{2}\left(\frac{b(L,u_{0})}{\lambda_{0}}\right)=1.$$
Therefore
$$\lambda_{0}=b(K+_{\phi}\varepsilon\cdot L,u_{0}).$$
Namely
$$b(K+_{\phi}\varepsilon\cdot L,u_{i})\rightarrow b(K+_{\phi}\varepsilon\cdot L,u_{0}).$$
as $i\rightarrow \infty.$

This shows that $K+_{\phi}\varepsilon\cdot L\in {\cal S}^{n}.$
\hfill$\Box$

{\bf Lemma 4.2}~ {\it If $K,L\in {\cal K}^{n}$, $\varepsilon>0$
and $\phi\in \Phi$, then
$$K+_{\phi}\varepsilon\cdot L\rightarrow K\eqno(4.1)$$ as} $\varepsilon\rightarrow
0^{+}$.

{\it Proof}~ From (3.4) and noting that $\phi_{2}$,
$\phi_{1}^{-1}$ and $b$ are continuous functions, we obtain
$$\lim_{\varepsilon\rightarrow 0^{+}}b(K+_{\phi}\varepsilon\cdot
L,u)=\lim_{\varepsilon\rightarrow
0^{+}}\frac{b(K,u)}{\displaystyle\phi_{1}^{-1}\left(1-\varepsilon\phi_{2}
\left(\frac{b(L,u)}{b(K+_{\phi}\varepsilon\cdot
L,u)}\right)\right)}.$$
Since $\phi^{-1}_{1}$ is continuous,
$\phi_{2}$ is bounded and in view of $\phi^{-1}_{1}(1)=1$, we have
$$\lim_{\varepsilon\rightarrow 0^{+}}\phi_{1}^{-1}\left(1-\varepsilon\phi_{2}
\left(\frac{b(L,u)}{b(K+_{\phi}\varepsilon\cdot
L,u)}\right)\right)=1.\eqno(4.2)$$ This yields
$$b(K+_{\phi}\varepsilon\cdot
L,u)\rightarrow b(K,u)$$ as $\varepsilon\rightarrow 0^{+}.$
\hfill$\Box$

{\bf Lemma 4.3}~ {\it If $K,L\in {\cal K}^{n}$, $0\leq i<n$ and
$\phi_{1}, \phi_{2}\in \Phi$, then}
$$\frac{d}{d\varepsilon}\bigg|_{\varepsilon=0^{+}}b(K+_{\phi}\varepsilon\cdot
L,u)^{n-i}=\frac{n-i}{(\phi_{1})'_{r}(1)}\cdot\phi_{2}
\left(\frac{b(L,u)}{b(K,u)}\right)\cdot b(K,u)^{n-i}.\eqno(4.3)$$

{\it Proof}~ Form (3.4), (4.1), Lemma 4.2 and notice that
$\phi_{1}^{-1}$, $\phi_{2}$ are continuous functions, we obtain
for $0\leq i<n$
$$\lim_{\varepsilon\rightarrow 0^+}\frac{b(K+_{\phi}\varepsilon\cdot
L,u)^{n-i}-b(K,u)^{n-i}}{\varepsilon}~~~~~~~~~~~~~~~~~~~~~~~~~~~~~~~~~~~~~~~~~~~~~~~~~~~~~~~~~~~~~~~~~~~$$
$$=\lim_{\varepsilon\rightarrow 0^+}\frac{\left(\frac{\displaystyle b(K,u)}{\displaystyle\phi_{1}^{-1}\left(1-\varepsilon\phi_{2}
\left(\frac{b(L,u)}{b(K+_{\phi}\varepsilon\cdot
L,u)}\right)\right)}\right)^{n-i}-b(K,u)^{n-i}}{\varepsilon}~~~~~~~~~~~~~~~~~~~$$
$$~~~~~=\lim_{\varepsilon\rightarrow 0^+}(n-i)\displaystyle b(K,u)^{n-i-1}\left(b(K,u)\phi_{2}\left(\frac{\displaystyle b(L,u)}
{\displaystyle b(K+_{\phi}\varepsilon\cdot
L,u)}\right)\right)\lim_{y\rightarrow
1^+}\frac{\phi_{1}^{-1}(y)-\phi_{1}^{-1}(1)}{y-1}$$
$$=\frac{n-i}{(\phi_{1})'_{r}(1)}\cdot\phi_{2}
\left(\frac{b(L,u)}{b(K,u)}\right)\cdot
b(K,u)^{n-i},~~~~~~~~~~~~~~~~~~~~~~~~~~~~~~~~~~~~~~~~~~~~~~~~$$
where
$$y=1-\varepsilon\phi_{2}
\left(\frac{b(L,u)}{b(K+_{\phi}\varepsilon\cdot L,u)}\right),$$
and note that $y\rightarrow 1^{+}$ as $\varepsilon\rightarrow
0^{+}.$
\hfill$\Box$

{\bf Lemma 4.4}~ {\it If $\phi\in \Phi_{2}$, $0\leq i<n$ and
$K,L\in {\cal K}^{n}$, then}
$$\frac{(\phi_{1})'_{r}(1)}{n-i}\cdot\frac{d}{d\varepsilon}\bigg|_{\varepsilon=0^{+}}A_{i}
(K+_{\phi}\varepsilon\cdot L) =\frac{1}{n}\int_{S^{n-1}}\phi_{2}
\left(\frac{b(L,u)}{b(K,u)}\right)\cdot b(K,u)^{n-i}
dS(u).\eqno(4.4)$$

{\it Proof}~ This follows immediately from (2.1) and Lemma
4.2.\hfill$\Box$

Denoting by $A_{\phi,i}(K,L)$, for any $\phi\in\Phi$ and $0\leq
i<n$, the integral on the right-hand side of (4.4) with $\phi_{2}$
replaced by $\phi$, we see that either side of the equation (4.4)
is equal to $A_{\phi_{2},i}(K,L)$ and hence this new Orlicz mixed
width integrals $A_{\phi,i}(K,L)$ has been born.

{\bf Definition 4.5}~ For $\phi\in \Phi$ and $0\leq i<n$, Orlicz
mixed width integrals of convex bodies $K$ and $L$,
$A_{\phi,i}(K,L)$, defined by
$$A_{\phi,i}(K,L)=:\frac{1}{n}\int_{S^{n-1}}\phi
\left(\frac{b(L,u)}{b(K,u)}\right)\cdot b(K,u)^{n-i}dS(u).
\eqno(4.5)$$

{\bf Lemma 4.6}~ {\it If $\phi_{1}, \phi_{2}\in \Phi$, $0\leq i<n$
and $K,L\in {\cal K}^{n}$, then}
$$A_{\phi_{2},i}(K,L)=\frac{(\phi_{1})'_{r}(1)}{n-i}\lim_{\varepsilon\rightarrow 0^+}
\frac{A_{i}(K+_{\phi}\varepsilon\cdot
L)-A_{i}(K)}{\varepsilon}.\eqno(4.6)$$

{\it Proof}~ This follows immediately from Lemma 4.4 and
(4.5).\hfill$\Box$

{\bf Lemma 4.7} {\it If $K,L\in {\cal K}^{n}$, $\phi\in \Phi$ and
any $A\in{\rm SL(n)}$, then for $\varepsilon>0$}
$$A(K+_{\phi}\varepsilon\cdot L)=(AK)+_{\phi}\varepsilon\cdot(AL).\eqno(4.7)$$

{\it Proof}~ For any $A\in{\rm SL(n)}$, from (2.1), we have
$$b(AK,u)=b(K,A^{t}u),~~ b(AL,u)=\rho(L,A^{t}u),~~ b(A(K\widehat{+}_{\phi}\varepsilon\cdot L,u)=b((K\widehat{+}_{\phi}\varepsilon\cdot L),A^{t}u).$$
Hence
$$b((AK+_{\phi}\varepsilon\cdot AL),u)~~~~~~~~~~~~~~~~~~~~~~~~~~~~~~~~~~~~~~~~~~~~~~~~~~~~~~~~$$
$$=\sup\left\{\lambda>0:\phi\left(\frac{b(AK,u)}{\lambda}\right)+\phi\left(\frac{b(AL,u)}{\lambda}\right)\leq
1\right\}$$
$$~=\sup\left\{\lambda>0:\phi\left(\frac{b(K,A^{t}u)}{\lambda}\right)+\phi\left(\frac{b(L,A^{t}u)}{\lambda}\right)\leq
1\right\}$$
$$=b(K+_{\phi}\varepsilon\cdot L,A^{t}u)~~~~~~~~~~~~~~~~~~~~~~~~~~~~~~~~~~~~~~~~~$$
$$~~~~~~~~~~~~~~~~~~~~~~~~~~~~~~~~~~~~~~~~~~~~~~~~~~=b(A(K+_{\phi}\varepsilon\cdot L),u).~~~~~~~~~~~~~~~~~~~~~~~~~~~~~~~~~~~~~~~~~~~~~~~~~~~~~~~~~~~~~~~~~~~~~~~~~~~~~~~~~~~~~~~~~~~~~~~~~~~~~~~~~~~~~~\Box$$

{\bf Lemma 4.8}~ {\it If $\phi\in \Phi$, $0\leq i<n$ and $K,L\in
{\cal K}^{n}$, then for $A\in SL(n)$,}
$$A_{\phi,i}(AK,AL)=A_{\phi,i}(K,L).\eqno(4.8)$$

{\it Proof}~ From (1.7), Lemma 3.3 and Lemma 4.7, we have for
$A\in{\rm SL(n)}$,
$$A_{\phi,i}(AK,AL)=\frac{(\phi_{1})
'_{r}(1)}{n-i}\cdot\frac{d}
{d\varepsilon}\bigg|_{\varepsilon=0^{+}}A_{i} (AK+
_{\phi}\varepsilon\cdot AL)~~~~~~~~~~~~~~~~$$
$$~~=\frac{(\phi_{1})
'_{r}(1)}{n-i}\cdot\frac{d}
{d\varepsilon}\bigg|_{\varepsilon=0^{+}}A_{i}
(A(K+_{\phi}\varepsilon\cdot L))$$
$$=\frac{(\phi_{1})
'_{r}(1)}{n-i}\cdot\frac{d}
{d\varepsilon}\bigg|_{\varepsilon=0^{+}}A_{i}
(K+_{\phi}\varepsilon\cdot L)~~$$
$$~~~~~~~~~~~~~~~~~~~~~~~~~~~~~~~~~~~~~~~~~~~=A_{\phi,i
}(K,L).~~~~~~~~~~~~~~~~~~~~~~~~~~~~~~~~~~~~~~~~~~~~~~~~~~~~~~~~~~~~~~~~~~~~~~~~\Box$$

\vskip 10pt \noindent{\large \bf 5 ~Orlicz width Minkowski
inequality}\vskip 10pt

In this section, we need define a Borel measure in $S^{n-1}$,
denoted by $A_{n,i}(K,\upsilon),$ call as width measure of convex
body $K$.

{\bf Definition 5.1}~ Let $K\in {\cal K}^{n}$ and $0\leq i<n$, the
width measure, denoted by $A_{n,i}(K,\upsilon),$ defined by
$$dA_{n,i}(K,\upsilon)=\frac{b(K,\upsilon)^{n-i}}{nA_{i}(K)}dS(\upsilon).
\eqno(5.1)$$

{\bf Lemma 5.2}~ (Jensen's inequality) {\it Let $\mu$ be a
probability measure on a space $X$ and $g: X\rightarrow I\subset
{\Bbb R}$ is a $\mu$-integrable function, where $I$ is a possibly
infinite interval. If $\psi: I\rightarrow {\Bbb R}$ is a convex
function, then
$$\int_{X}\psi(g(x))d\mu(x)\geq\psi\left(\int_{X}g(x)d\mu(x)\right).
\eqno(5.2)$$ If $\psi$ is strictly convex, equality holds if and
only if $g(x)$ is constant for $\mu$-almost all $x\in X$} (see
[12, p.165]).

{\bf Lemma 5.3}~ {\it Suppose that $\phi: [0,\infty)\rightarrow
(0,\infty)$ is decreasing and convex with $\phi(0)=\infty$. If
$K,L\in{\cal K}^{n}$ and $0\leq i<n$, then
$$\frac{1}{nA_{i}(K)}\int_{S^{n-1}}\phi
\left(\frac{b(L,u)}{b(K,u)}\right)b(K,u)^{n-i}dS(u)\geq
\phi\left(\left(\frac{A_{i}(L)}{A_{i}(K)}\right)^{1/(n-i)}\right)
.\eqno(5.3)$$ If $\phi$ is strictly convex, equality holds if and
only if $K$ and $L$ have similar width.}

{\it Proof}~ For $K\in {\cal K}^{n-1}$, $0\leq i<n$ and any $u\in
S^{n-1}$, since
$$\int_{S^{n-1}}dA_{n,i}(K,\upsilon)=1,$$
so the width measure $\frac{\displaystyle
b(K,u)^{n-i}}{\displaystyle nA_{i}(K)}dS(u)$ is a probability
measure on $S^{n-1}$.

Hence, from (2.4), (2.5), (5.1) and by using Jensen's inequality
(5.2), and in view of $\phi$ is decreasing, we obtain
$$\frac{1}{nA_{i}(K)}\int_{S^{n-1}}\phi
\left(\frac{b(L,u)}{b(K,u)}\right)b(K,u)^{n-i}dS(u)$$
$$~~=\int_{S^{n-1}}
\phi\left(\frac{b(L,u)}{b(K,u)}\right)dA_{n,i}(K,u)$$
$$\geq\phi\left(\frac{A_{i}(K,L)}{A_{i}(K)}\right)~~~~~~~~~~~~~~~~~$$
$$\geq\phi\left(\left(\frac{A_{i}(L)}{A_{i}(K)}\right)^{1/(n-i)}\right).~~~~~~$$
Next, we discuss the equal condition of (5.3). Suppose the
equality holds in (5.3) and $\phi$ is strictly convex, form the
equality condition of (2.5), so there exist $r>0$ such that
$$b(L,u)=rb(K,u),$$ for all $u\in S^{n-1}$.
On the other hand, suppose that $K$ and $L$ have similar width,
i.e. there exist $\lambda>0$ such that $b(L,u)=\lambda b(K,u)$ for
all $u\in S^{n-1}$. Hence
$$\frac{1}{nA_{i}(K)}\int_{S^{n-1}}\phi
\left(\frac{b(L,u)}{b(K,u)}\right)b(K,u)^{n-i}dS(u)$$
$$~~~~~~~~~~~~~~~~~~~~~~~~~~~~~~~~~~~~~~~~~=\frac{1}{nA_{i}(K)}\int_{S^{n-1}}\phi
\left(\left(\frac{A_{i}(L)}{A_{i}(K)}\right)^{1/(n-i)}\right)b(K,u)^{n-i}dS(u)$$
$$~=\phi\left(\left(\frac{A_{i}(L)}{A_{i}(K)}\right)^{1/(n-i)}\right).$$
This implies the equality in (5.3) holds.
\hfill$\Box$

{\bf Theorem 5.4}~ (Orlicz width Minkowski inequality) {\it If
$K,L\in {\cal K}^{n}$, $0\leq i<n$ and $\phi\in \Phi$, then
$$A_{\phi,i}(K,L)\geq
A_{i}(K)\cdot\phi\left(\left(\frac{A_{i}(L)}{A_{i}(K)}\right)^{1/(n-i)}\right)
.\eqno(5.4)$$ If $\phi$ is strictly convex, equality holds if and
only if $K$ and $L$ have similar width.}

{\it Proof}~ This follows immediately from (4.5) and Lemma 5.3.
\hfill$\Box$

{\bf Corollary 5.5} {\it If $K,L\in {\cal K}^{n}$, $0\leq i<n$ and
$p\geq 1$, then
$$A_{-p,i}(K,L)^{n-i}\geq A_{i}(K)^{n-i+p}A_{i}(L)^{-p},\eqno(5.5)$$
with equality if and only if $K$ and $L$ have similar width.}

{\it Proof}~ This follows immediately from Theorem 5.4 with
$\phi_{1}(t)=\phi_{2}(t)=t^{-p}$ and $p\geq 1$.
\hfill$\Box$

Taking $i=0$ in (5.6), this yields $L_{p}$-Minkowski inequality is
following: If $K,L\in {\cal K}^{n}$ and $p\geq 1$, then
$$A_{-p}(K,L)^{n}\geq A(K)^{n+p}\cdot A(L)^{-p},$$
with equality if and only if $K$ and $L$ have similar width.

{\bf Corollary 5.6} {\it Let $K,L\in {\cal M}\subset{\cal K}^{n}$,
$0\leq i<n$ and $\phi\in \Phi$, and if either
$$A_{\phi,i}(Q,K)=A_{\phi,i}(Q,L),~ {\it for~ all}~ Q\in{\cal M}\eqno(5.6)$$
or
$$\frac{A_{\phi,i}(K,Q)}{A_{i}(K)}=\frac{A_{\phi,i}(L,Q)}{A_{i}(L)},~ {\it for~ all}~ Q\in{\cal M},\eqno(5.7)$$
then} $K=L.$

{\it Proof}~ Suppose (5.6) hold. Taking $K$ for $Q$, then from
(2.3), (4.5) and (5.3), we obtain
$$A_{i}(K)=A_{\phi,i}(K,L)\geq A_{i}(K)\phi\left(\left(\frac{A_{i}(L)}{A_{i}(K)}\right)^{1/(n-i)}\right)$$
with equality if and only if $K$ and $L$ have similar width. Hence
$$1\geq\phi\left(\left(\frac{A_{i}(L)}{A_{i}(K)}\right)^{1/(n-i)}\right),$$
with equality if and only if $K$ and $L$ have similar width. Since
$\phi$ is decreasing function on $(0,\infty),$ this follows that
$$A_{i}(K)\leq A_{i}(L),$$
with equality if and only if $K$ and $L$ have similar width. On
the other hand, if taking $L$ for $Q$, we similar get
$A_{i}(K)\geq A_{i}(L),$ with equality if and only if $K$ and $L$
have similar width. Hence $A_{i}(K)=A_{i}(L),$ and $K$ and $L$
have similar width, it follows that $K$ and $L$ must be equal.

Suppose (5.7) hold. Taking $L$ for $Q$, then from from (2.3),
(4.5) and (5.3), we obtain
$$1=\frac{A_{\phi,i}(K,L)}{A_{i}(K)}\geq
\phi\left(\left(\frac{A_{i}(L)}{A_{i}(K)}\right)^{1/(n-i)}\right),$$
with equality if and only if $K$ and $L$ are dilates. Hence
$$1\geq\phi\left(\left(\frac{A_{i}(L)}{A_{i}(K)}\right)^{1/(n-i)}\right),$$
with equality if and only if $K$ and $L$ have similar width. Since
$\phi$ is decreasing function on $(0,\infty),$ this follows that
$$A_{i}(K)\leq A_{i}(L),$$
with equality if and only if $K$ and $L$ have similar width. On
the other hand, if taking $K$ for $Q$, we similar get
$A_{i}(K)\geq A_{i}(L),$ with equality if and only if $K$ and $L$
have similar width. Hence $A_{i}(K)=A_{i}(L),$ and $K$ and $L$
have similar width, it follows that $K$ and $L$ must be equal.
\hfill$\Box$

When $\phi_{1}(t)=\phi_{2}(t)=t^{-p}$ and $p\geq 1$, Corollary 5.6
becomes the following result.

{\bf Corollary 5.7} {\it Let $K,L\in {\cal M}\subset{\cal K}^{n}$,
$0\leq i<n$ and $p\geq 1$, and if either
$$A_{-p,i}(K,Q)=A_{-p,i}(L,Q),~ {\it for~ all}~ Q\in{\cal M}$$
or
$$\frac{A_{-p,i}(K,Q)}{A_{i}(K)}=\frac{A_{-p,i}(L,Q)}{A_{i}(L)},~ {\it for~ all}~ Q\in{\cal M},$$
then} $K=L.$

\vskip 10pt \noindent{\large \bf 6 ~Orlicz width Brunn-Minkowski
inequality}\vskip 10pt

{\bf Lemma 6.1}~ {\it If $K,L\in {\cal K}^{n}$, $0\leq i<n$, and
$\phi_{1}, \phi_{2}\in \Phi$, then}
$$A_{i}(K+_{\phi}L)=A_{\phi_{1},i}(K+_{\phi}L, K)
+A_{\phi_{2},i}(K+_{\phi}L, L).\eqno(6.1)$$

{\it Proof}~ From (2.3), (3.1), (3.4) and (4.5), we have for
$K+_{\phi}L=Q\in {\cal K}^{n}$
$$A_{\phi_{1},i}(Q,K)+A_{\phi_{2},i}(Q,L)~~~~~~~~~~~~~~~~~~~~~~~~~~~~~~~~~~~~~~~~~~~~~~~~~~~~~~~~~~~~~~~~~~~~$$
$$=\frac{1}{n}\int_{S^{n-1}}\left(\phi_{1}
\left(\frac{b(K,u)}{b(Q,u)}\right)+\phi_{2}
\left(\frac{b(L,u)}{b(Q,u)}\right)\right)b(Q,u)^{n-i}dS(u)$$
$$=\frac{1}{n}\int_{S^{n-1}}\phi\left(\frac{b(K,u)}{b(Q,u)},\frac{b(L,u)}{b(Q,u)}\right)
b(Q,u)^{n-i}dS(u)~~~~~~~~~~~~~~~~$$
$$=A_{i}(Q).~~~~~~~~~~~~~~~~~~~~~~~~~~~~~~~~~~~~~~~~~~~~~~~~~~~~~~~~~~~~~~~~~~~~\eqno(6.2)$$
This completes the proof.\hfill$\Box$

{\bf Theorem 6.2}~ (Orlicz width Brunn-Minkowski inequality)~ {\it
If $K,L\in{\cal K}^{n}$, $0\leq i<n$ and $\phi\in \Phi_{2}$, then
$$1\geq\phi\left(\left(\frac{A_{i}(K)}
{A_{i}(K+_{\phi}L)}\right)^{1/(n-i)},\left(\frac{A_{i}(L)}{A_{i}(K{+}_
{\phi}L)}\right)^{1/(n-i)}\right).\eqno(6.3)$$ If $\phi$ is
strictly convex, equality holds if and only if $K$ and $L$ have
similar width.}

{\it Proof}~ From (5.4) and Lemma 6.1, we have
$$A_{i}(K+_{\phi}L)=A_{\phi_{1},i}(K+_{\phi}L, K)
+A_{\phi_{2},i}(K+_{\phi}L,
L)~~~~~~~~~~~~~~~~~~~~~~~~~~~~~~~~~~~~~~~~~~~~~~~~~$$
$$~~~~~~~~~~~~~~~~~~\geq A_{i}(K+_{\phi}L)\left(\phi_{1}\left(
\left(\frac{A_{i}(K)}{A_{i}(K+_{\phi}L)}\right)^{1/(n-i)}\right)+
\phi_{2}\left(\left(\frac{A_{i}(L)}{A_{i}(K+_{\phi}L)}\right)^{1/(n-i)}\right)\right)$$
$$~=A_{i}(K+_{\phi}L)\phi\left(\left(\frac{A_{i}(K)}
{A_{i}(K+_{\phi}L)}\right)^{1/(n-i)},\left(\frac{A_{i}(L)}{A_{i}(K{+}_
{\phi}L)}\right)^{1/(n-i)}\right).$$ This is just inequality
(6.3). From the equality condition of (5.4), if follows that if
$\phi$ is strictly convex, equality in (6.3) holds if and only if
$K$ and $L$ have similar width.
\hfill$\Box$

{\bf Corollary 6.3} {\it If $K,L\in {\cal K}^{n}$, $0\leq i<n$ and
$p\geq 1$, then
$$A_{i}(K+_{p}L)^{-p/(n-i)}\geq A_{i}(K)^{-p/(n-i)}+A_{i}(L)^{-p/(n-i)},\eqno(6.4)$$
with equality if and only if $K$ and $L$ have similar width.}

{\it Proof} The result follows immediately from Theorem 6.2 with
$\phi(x_{1},x_{2})=x_{1}^{-p}+x_{2}^{-p}$ and $p\geq 1$.
\hfill$\Box$

Taking $i=0$ in (6.4), this yields the $L_{p}$-Brunn-Minkowski
inequality for the width integrals. If $K,L\in{\cal K}^{n}$ and
$p\geq 1$, then
$$A(K+_{p}L)^{-p/n}\geq A(K)^{-p/n}+A(L)^{-p/n},$$
with equality if and only if $K$ and $L$ have similar width.

{\bf Corollary 6.4}~ {\it If $K,L\in {\cal K}^{n}$, $0\leq i<n$
and $\phi\in \Phi$, then
$$A_{\phi,i}(K,L)\geq
A_{i}(K)\cdot\phi\left(\left(\frac{A_{i}(L)}{A_{i}(K)}\right)^{1/(n-i)}\right)
.\eqno(6.5)$$ If $\phi$ is strictly convex, equality holds if and
only if $K$ and $L$ have similar width.}

{\it Proof}~ Let
$$K_{\varepsilon}=K+_{\phi}\varepsilon\cdot L.$$
From (4.6) and in view of the Orlicz-Brunn-Minkowski inequality
(6.3), we obtain
$$\frac{n-i}{(\phi_{1})'_{r}(1)}\cdot A_{\phi_{2},i}(K,L)=\frac{d}{d\varepsilon}\bigg|_{\varepsilon=0^{+}}A_{i}
(K_{\varepsilon})~~~~~~~~~~~~~~~~~~~~~~~~~~~~~~~~~~~~~~~~~~~~~~~~~~~~~~~~~~~~~~~~~~~~~~~~~~~~~~~~~~~~~~~~~~~~~~~$$
$$=\lim_{\varepsilon\rightarrow 0^{+}}\frac{A_{i}(K_{\varepsilon})-A_{i}(K)}{\varepsilon}~~~~~~~~~~~~~~~~~~~~~~~~~~~~~~~~~~~~~~~~~~~~~~~~$$
$$~~~~~~~~~~=\lim_{\varepsilon\rightarrow 0^{+}}\frac{\displaystyle 1-\frac{A_{i}(K)}{A_{i}(K_{\varepsilon})}}
{\displaystyle
1-\phi_{1}\left(\left(\frac{A_{i}(K)}{A_{i}(K_{\varepsilon})}\right)^{\frac{1}{n-i}}\right)}\cdot\frac{\displaystyle
1-\phi_{1}\left(\left(\frac{A_{i}(K)}{A_{i}(K_{\varepsilon})}\right)^{\frac{1}{n-i}}\right)}{\varepsilon}\cdot
A_{i}(K_{\varepsilon})$$
$$~~~~~~~~~~~~~~~~=\lim_{t\rightarrow 1^{-
}}\frac{1-t}{\displaystyle\phi_{1}(1)-\phi_{1}\left(t^{\frac{1}{n-i}}\right)}\cdot\lim_{\varepsilon\rightarrow
0^{+}}\frac{\displaystyle
1-\phi_{1}\left(\left(\frac{A_{i}(K)}{A_{i}(K_{\varepsilon})}\right)^{\frac{1}{n-i}}\right)}{\varepsilon}
\cdot\lim_{\varepsilon\rightarrow 0^{+}}A_{i}(K_{\varepsilon})~~$$
$$\geq\frac{n-i}{(\phi_{1})'_{r}(1)}\cdot\lim_{\varepsilon\rightarrow
0^{+}}\phi_{2}\left(\left(\frac{A_{i}(L)}{A_{i}(K_{\varepsilon})}\right)^{\frac{1}{n-i}}\right)
\cdot\lim_{\varepsilon\rightarrow
0^{+}}A_{i}(K_{\varepsilon})~~~~~~~~~~$$
$$=\frac{n-i}{(\phi_{1})'_{r}(1)}\cdot\phi_{2}\left(\left(\frac{A_{i}(L)}{A_{i}(K)}\right)^{\frac{1}{n-i}}\right)
\cdot A_{i}(K).~~~~~~~~~~~~~~~~~~~~~~~~~\eqno(6.6)$$ Replacing
$\phi_{2}$ by $\phi$, (6.6) becomes (6.5). If $\phi$ is strictly
convex, from the equality condition of (6.3), it follows that the
equality holds in (6.5) if and only if $K$ and $L$ have similar
width.

This proof is
complete.\hfill$\Box$

{\bf Acknowledgements}

The author thanks Ms. Jinhua Ji for valuable comments and suggestions which improved the presentation of this
manuscript.


\begin{thebibliography}{zz}
\vskip 6pt {\small

\bibitem{35-1} A. D. Aleksandrov, Zur Theorie der gemischten Volumina von konvexen
K\"{o}rpern, I: Verall-gemeinerung einiger Begriffe der Theorie
der konvexen K\"{o}rper, {\it Mat. Sbornik  N. S.} {\bf 2}, (1937)
947-972.
\bibitem{36-2} W. Fenchel, B. Jessen, Mengenfunktionen und konvexe K\"{o}rper,
{\it Danske Vid Selskab Mat-fys Medd}, {\bf 16} (1938), 1-31.
\bibitem{1-3} W. J. Firey, Polar means of convex bodies and a dual to the
Brunn-Minkowski theorem, {\it Canad. J. Math.,} {\bf 13} (1961),
444-453.
\bibitem{2-4} W. J. Firey, p-means of convex bodies, {\it Math.
Scand.,} {\bf 10} (1962), 17-24.
\bibitem{38-5} R. J. Gardner, Geometric Tomography, Cambridge University Press, second edition, New York, 2006.
\bibitem{21-6} R. J. Gardner, D. Hug, W. Weil, The Orlicz-Brunn-Minkowski theory: a general framework, additions, and
inequalities, {\it J. Diff. Geom.,} {\bf 97}(3) (2014), 427-476.
\bibitem{24-7} C. Haberl, E. Lutwak, D. Yang, G. Zhang, The even Orlicz Minkowski problem, {\it Adv. Math.,} {\bf 224} (2010),
2485-2510.
\bibitem{3-8} C. Haberl, L. Parapatits, The Centro-Affine Hadwiger Theorem, {\it J. Amer. Math.
Soc.,} in press.
\bibitem{4-9} C. Haberl, F. E. Schuster, Asymmetric affine $L_{p}$ Sobolev inequalities,
{\it J. Funct. Anal.,} {\bf 257} (2009), 641-658.
\bibitem{5-10} C. Haberl, F. E. Schuster, General $L_{p}$ affine isoperimetric
inequalities, {\it J. Differential Geom.,} {\bf 83} (2009), 1-26.
\bibitem{6-11} C. Haberl, F. E. Schuster, J. Xiao, An asymmetric affine P\'{o}lya-Szeg\"{o}
principle, {\it Math. Ann.,} {\bf 352} (2012), 517-542.
\bibitem{39-12} J. Hoffmann-J$\phi$gensen, Probability With a View Toward
Statistics, Vol. I, Chapman and Hall, New York, 1994, 165-243.
\bibitem{25-13} Q. Huang, B. He, On the Orlicz Minkowski problem for polytopes, {\it Discrete Comput. Geom.,} {\bf 48} (2012), 281-297.
\bibitem{26-14} M. A. Krasnosel'skii, Y. B. Rutickii, Convex Functions and Orlicz Spaces, P. Noordhoff Ltd., Groningen, 1961.
\bibitem{7-15} M. Ludwig, M. Reitzner, A classification of $SL(n)$ invariant valuations, {\it Ann. Math.,} {\bf 172} (2010),
1223-1271.
\bibitem{8-16} E. Lutwak, The Brunn-Minkowski-Firey theory I. mixed volumes and the Minkowski problem. {\it J. Diff. Goem.,} {\bf 38} (1993), 131-150.
\bibitem{34-17} E. Lutwak, Mixed width-integrals of convex bodies, {\it Israel J. Math.,} {\bf 28}
(1977), 249-253.
\bibitem{9-18} E. Lutwak, D. Yang, G. Zhang, On the $L_{p}$-Minkowski problem, {\it Trans.
Amer. Math. Soc.,} {\bf 356} (2004), 4359-4370.
\bibitem{10-19} E. Lutwak, D. Yang, G. Zhang, $L_{p}$ John ellipsoids, {\it Proc. London Math. Soc.,} {\bf 90} (2005),
497-520.
\bibitem{11-20} E. Lutwak, D. Yang, G. Zhang, $L_{p}$ affine isoperimetric
inequalities, {\it J. Differential Geom.,} {\bf 56} (2000),
111-132.
\bibitem{12-21} E. Lutwak, D. Yang, G. Zhang, Sharp affine $L_{p}$ Sobolev
inequalities, {\it J. Differential Geom.,} {\bf 62} (2002), 17-38.
\bibitem{13-22} E. Lutwak, D. Yang, G. Zhang, The Brunn-Minkowski-Firey
inequality for nonconvex sets, {\it Adv. Appl. Math.,} {\bf 48}
(2012), 407-413.
\bibitem{19-23} E. Lutwak, D. Yang, G. Zhang, Orlicz projection bodies, {\it Adv. Math.,} {\bf 223} (2010), 220-242.
\bibitem{20-24} E. Lutwak, D. Yang, G. Zhang, Orlicz centroid bodies, {\it J. Differential Geom.,} {\bf 84} (2010), 365-387.
\bibitem{14-25} L. Parapatits, $SL(n)$-Covariant $L_{p}$-Minkowski Valuations, {\it J. Lond. Math. Soc.,} in press.
\bibitem{15-26} L. Parapatits, $SL(n)$-Contravariant $L_{p}$-Minkowski Valuations, {\it Trans. Amer. Math. Soc.}, in press.
\bibitem{27-27} M. M. Rao and Z. D. Ren, Theory of Orlicz Spaces, Marcel Dekker, New York, 1991.
\bibitem{37-28} R. Schneider, Convex Bodies: The Brunn-Minkowski Theory, Second Edition, Cambridge Univ. Press, 2014.
\bibitem{16-29} C. Sch\"{u}tt, E. Werner, Surface bodies and $p$-affine surface
area, {\it Adv. Math.,} {\bf 187} (2004), 98-145.
\bibitem{17-30} E. M. Werner, R\'{e}nyi divergence and $L_{p}$-affine surface area for
convex bodies, {\it Adv. Math.,} {\bf 230} (2012), 1040-1059.
\bibitem{18-31} E. Werner, D. P. Ye, New $L_{p}$ affine isoperimetric inequalities,
{\it Adv. Math.,} {\bf 218} (2008), 762-780.
\bibitem{22-32} D. Xi, H. Jin, G. Leng, The Orlicz Brunn-Minkwski inequality, {\it Adv.
Math.,} {\bf 260} (2014), 350-374.
\bibitem{28-33} D. Ye, Dual Orlicz-Brunn-Minkowski theory: dual Orlicz
$L_\varphi$ affine and geominimal surface areas, {\it J. Math.
Anal. Appl.,} {\bf 443} (2016), 352-371.
\bibitem{30-34} C.-J. Zhao, On the Orlicz-Brunn-Minkowski theory, {\it Balkan J.
Geom. Appl.}, {\bf 22} (2017), 98-121.
\bibitem{31-35} C.-J. Zhao, Orlicz dual mixed volumes, {\it Results Math.}, {bf
68} (2015), 93-104.
\bibitem{32-36} C.-J. Zhao, Orlicz dual affine quermassintegrals, {\it Forum Math.}, 2018, DOI: https://doi.org/10.1515/forum-2017-0174.
\bibitem{33-37} C.-J. Zhao, Orlicz-Brunn-Minkowski inequality for radial Blaschke-Minkowski homomorphisms, {\it Quaestiones Math.}, 2018, in press.
\bibitem{29-38} B. Zhu, J. Zhou, W. Xu, Dual Orlicz-Brunn-Minkwski theory, {\it Adv.
Math.,} {\bf 264} (2014), 700-725.
\bibitem{23-39} G. Zhu, The Orlicz centroid inequality for convex bodies, {\it Adv.
Appl. Math.,} {\bf 48} (2012), 432-445.}

\end{thebibliography}
\end{document}